\documentclass[english]{article}
\usepackage[T1]{fontenc}
\usepackage[latin9]{inputenc}
\usepackage{geometry}
\usepackage{amsmath}

\makeatletter
\date{}

\makeatother

\usepackage{babel}
\begin{document}

\title{Discriminants of Chebyshev-like polynomials and their generating
functions}

\author{Khang Tran}
\maketitle
\begin{abstract}
In his paper of 2000, Kenneth B. Stolarsky made various observations
and conjectures about discriminants and generating functions of certain
types of Chebyshev-like polynomials. We prove several of these conjectures.
One of our proofs involves Wilf-Zeilberger pairs and a contiguous
relation for hypergeometric series.
\end{abstract}

\section{Introduction}

We begin by recalling the notion of discriminant. Suppose that $P_{m}(x)$
is a polynomial of degree $m$ whose roots are $x_{1}$, $x_{2}$,
$\ldots$, $x_{m}$. Then the disciminant of $P_{m}(x)$, which we
will denote by $\Delta_{x}P_{m}(x)$, is
\[
\Delta_{x}P_{m}(x)=a_{m}^{2m-2}\prod_{1\le i<j\le m}(x_{i}-x_{j})^{2}
\]
where $a_{m}$ is the leading coefficient. If $P_{m}(x)$ is monic,
it is easy to prove the useful fact that

\begin{equation}
\Delta_{x}P_{m}(x)=(-1)^{m(m-1)/2}P'(x_{1})\cdots P'(x_{m}).
\end{equation}

The study of discriminants and resultants of specific types of polynomials
has a long history, and includes contributions from both Stieltjes
and Hilbert (see \cite{aar} and \cite{gkz} respectively). Formulas
for various specific types of polynomials are given, e.g., in \cite{aar},
\cite{apostol}, \cite{dilcherstolarsky}, \cite{gkz}, and \cite{gishe-ismail}.
The importance of this subject is discussed in the introduction to
\cite{gishe-ismail}. Here we contribute to this area by establishing
some of the conjectures in \cite{stolarsky}.

The discriminant of the product of the two polynomials 
\[
K_{n}(x,q)=(1+x)^{2n}+qx^{n}
\]
and 
\[
f_{m}(x)=(x^{2m+1}-1)/(x-1)
\]
 has some remarkable properties in terms of roots, generating functions
and divisibility. The following observation of Stolarsky limits the
possible range of the roots of the resulting polynomial \cite{stolarsky}.

\paragraph{Proposition. (K. Stolarsky, 2000).}

Let $a<b$ be real. Let $K(x,q)$ be any polynomial in $x$ and $q$
such that

\,\,\,(i) $K(x,q)$ is never zero for $|x|=1$ unless $a\leq q\leq b$;

\,\,\,(ii) $K(x,q)$ has no multiple zeros unless $a\leq q\leq b$.\\
Say $f_{m}(x)$ is a sequence of polynomials with no multiple zeros,
and the zeros of $f_{m}(x)$ all lie on $|x|=1$. Then
\[
g_{m}(q):=\Delta_{x}(K(x,q)f_{m}(x))
\]
 is a sequence of polynomials in $q$ whose roots are all in $[a,b]$. 

\medskip{}
An application of this theorem shows that the discriminant of the
product of the two polynomials $K_{n}(x,q)$ and $f_{m}(x)$ defined
above has all real roots in $[-2^{2n},0]$. Calculation of this discriminant
gives rise to the following polynomial:
\[
H_{m}^{(n)}(q):=\prod_{k=1}^{m}(q+(2\cos k\theta_{m}+2)^{n})
\]
where 
\[
\theta_{m}=2\pi/(2m+1).
\]

See Table 1 for a short tabulation of these polynomials for $1\le n\le3$.
It immediately suggests that for $n=1$ we have a familiar type of
Fibonacci polynomial (for $q=1$ the values are 1, 2, 5, 13, 34, etc.).
However, for $n\ge2$ the nature of these polynomials is much less
obvious.

\begin{table}
\begin{itemize}
\item \noindent $n=1$\\
\\
$H_{0}^{(1)}(q)=1$\\
$H_{1}^{(1)}(q)=1+q$\\
$H_{2}^{(1)}(q)=1+3q+q^{2}$\\
$H_{3}^{(1)}(q)=1+6q+5q^{2}+q^{3}$\\
$H_{4}^{(1)}(q)=1+10q+15q^{2}+7q^{3}+q^{4}$\\
\rule[0.5ex]{1\linewidth}{1pt}
\item \noindent $n=2$\\
\\
$H_{0}^{(2)}(q)=1$\\
$H_{1}^{(2)}(q)=1+q$\\
$H_{2}^{(2)}(q)=1+7q+q^{2}$\\
$H_{3}^{(2)}(q)=1+26q+13q^{2}+q^{3}$\\
$H_{4}^{(2)}(q)=1+70q+87q^{2}+19q^{3}+q^{4}$\\
\rule[0.5ex]{1\linewidth}{1pt}
\item $n=3$\\
\\
$H_{0}^{(3)}(q)=1$\\
$H_{1}^{(3)}(q)=1+q$\\
$H_{2}^{(3)}(q)=1+18q+q^{2}$\\
$H_{3}^{(3)}(q)=1+129q+38q^{2}+q^{3}$\\
$H_{4}^{(3)}(q)=1+571q+627q^{2}+58q^{3}+q^{4}$
\end{itemize}
\caption{The polynomials $H_{m}^{(n)}(q)$.}
\end{table}

This polynomial has several interesting properties conjectured by
Stolarsky. Some of these properties are clear from proposition 2.3.
The previously conjectured formulas for the generating functions of
$H_{m}^{(1)}(q)$ and $H_{m}^{(2)}(q)$ are obtained in Sections 3
and 4. It is still unknown to the author what the generating function
of $H_{m}^{(n)}(q)$ is for $n\ge3$. Stolarsky \cite{stolarsky}
conjectured the following generating function for $H_{m}^{(3)}(q)$:
\[
\frac{(1-t)^{7}-qt^{2}(1-t)(t+3)(1+3t)}{(1-t)^{8}-qt(1-t)^{2}(1+14t+34t^{2}+14t^{3}+t^{4})+x^{2}t^{4}}.
\]

To find these generating functions, we use some knowledge about hypergeometric
series and the Wilf-Zeilberger algorithm. A hypergeometric series
is defined as
\[
_{p}F_{q}\left(\begin{array}{c}
a_{1},a_{2},\ldots,a_{p}\\
b_{1},b_{2},\ldots,b_{q}
\end{array};x\right):=\sum_{n=0}^{\infty}\frac{(a_{1})_{n}(a_{2})_{n}\ldots(a_{p})_{n}}{(b_{1})_{n}(b_{2})_{n}\ldots(b_{q})_{n}}\frac{x^{n}}{n!}.
\]
Euler obtained the following contiguous relation {[}1, equation (2.5.3){]}
for a hypergeometric series $_{2}F_{1}$:
\begin{eqnarray*}
_{2}F_{1}\left(\begin{array}{c}
a,b\\
c
\end{array};u\right) & = & \frac{c+(a-b+1)u}{c}{}_{2}F_{1}\left(\begin{array}{c}
a+1,b\\
c+1
\end{array};u\right)\\
 &  & -\frac{(a+1)(c-b+1)u}{c(c+1)}{}_{2}F_{1}\left(\begin{array}{c}
a+2,b\\
c+2
\end{array};u\right).
\end{eqnarray*}
A Wilf-Zeilberger pair $(F,G)$ \cite{pwz} satisfies the equation
\[
F(m+1,i)-F(m,i)=G(m,i+1)-G(m,i).
\]
By telescoping summation, if $G(m,a)=G(m,b+1)=0$ then $\sum_{i=a}^{b}F(m,i)$
does not depend on $m$. We will use this information to prove the
identity
\[
\sum_{i=0}^{2k}(-1)^{i}\left(\begin{array}{c}
m+k+i\\
m+k-i
\end{array}\right)\left(\begin{array}{c}
m+3k-i\\
m-k+i
\end{array}\right)=(-1)^{k}\sum_{i=0}^{m}\left(\begin{array}{c}
2k\\
m-i
\end{array}\right)\left(\begin{array}{c}
4k+i\\
i
\end{array}\right)
\]
in section 5. This identity is crucial for finding the generating
function of $H_{m}^{(2)}(q)$.

\section{A general form of the discriminant}

In this section we will compute the disciminant of $K_{n}(x,q)$.
For further computations of the resultants and disciminants of different
kinds of polynomials, see \cite{apostol}, \cite{dilcherstolarsky},
\cite{gkz}, and \cite{roberts}. We first note that in the cases
$q=0$ and $q=-2^{2n}$ the polynomial $K_{n}(x,q)$ has multiple
roots. These roots arise solely from the $K_{n}(x,q)$ polynomial.
So one expects that $\Delta_{x}(K_{n}f_{m})$ has factors $q$ and
$(q+2^{2n})$ . Other factors of $\Delta_{x}(K_{n}f_{m})$ appear
when $K_{n}(x,q)$ and $f_{m}(x)$ have common roots.

\paragraph{Proposition 2.1.}

The discriminant of $K_{n}(x,q)$ is given by
\[
\Delta_{x}(K_{n}(x,q))=n^{2n}q^{2n-1}(q+2^{2n}).
\]
\\
Proof. Let $x_{l}$ be a root of $K_{n}(x)$ (we suppress the parameter
$q$ for a moment). Then
\begin{eqnarray*}
(1+x_{l})K_{n}'(x_{l}) & = & 2n(1+x_{l})^{2n}+nqx_{l}^{n-1}(1+x_{l})\\
 & = & -2nqx_{l}^{n}+nqx_{l}^{n-1}(1+x_{l})\\
 & = & nqx_{l}^{n-1}(1-x_{l}).
\end{eqnarray*}
Since $(-1)^{2n(2n-1)/2}=(-1)^{n}$, it follows that
\[
(-1)^{n}\Delta_{x}(K_{n}(x,q))\prod_{l=1}^{2n}(1+x_{l})=n^{2n}q^{2n}\prod_{l=1}^{2n}(1-x_{l}).
\]
Moreover we have
\[
\prod_{l=1}^{2n}(1+x_{l})=K_{n}(-1,q)=(-1)^{n}q
\]
and
\[
\prod_{l=1}^{2n}(1-x_{l})=K_{n}(1,q)=q+4^{n},
\]
and the proof follows.

\medskip{}

\paragraph{Proposition 2.2.}

The discriminant of $f_{m}(x)$ is
\[
\Delta_{x}f_{m}(x)=(-1)^{m}(2m+1)^{2m-1}.
\]
\\
Proof. According to (1) we have
\begin{eqnarray*}
\Delta_{x}f_{m}(x) & = & (-1)^{2m(2m-1)/2}\prod_{k=1}^{2m}f'(e^{ik\theta_{m}})\\
 & = & (-1)^{m}\prod_{k=1}^{2m}\frac{(2m+1)}{e^{ik\theta_{m}}-1}
\end{eqnarray*}
Since the denominator factors are the nonzero roots of $(x+1)^{2m+1}=1$,
their product is $2m+1$, and the proposition follows.

\paragraph*{Proposition 2.3. }

The discriminant of $K_{n}(x,q)f_{m}(x)$ is
\begin{eqnarray*}
\Delta_{x}(K_{n}f_{m}) & = & C_{m}^{(n)}q^{2n-1}(q+2^{2n})\prod_{k=1}^{m}(q+(2\cos k\theta_{m}+2)^{n})^{4}\\
 & = & C_{m}^{(n)}q^{2n-1}(q+2^{2n})H_{m}^{(n)}(q)^{4}
\end{eqnarray*}
where $ $ 
\[
C_{m}^{(n)}=(-1)^{m}(2m+1)^{2m-1}n^{2n}.
\]
 Proof. From the definition of discriminant we note that
\begin{eqnarray*}
\Delta_{x}(K_{n}f_{m}) & = & \Delta_{x}(K_{n})\Delta_{x}(f_{m})\prod_{\substack{1\le k\le2m\\
1\le l\le2n
}
}(x_{l}-e^{ik\theta_{m}})^{2}
\end{eqnarray*}
where the $x_{l}$'s are roots of $K_{n}$. Thus by proposition 2.1
and 2.2 it suffices to show that 
\begin{eqnarray*}
\prod_{l,1\le k\le2m}(x_{i}-e^{ik\theta_{m}})^{2} & = & \prod_{k=1}^{m}(q+(2\cos k\theta+2)^{n})^{4}\\
 & = & \prod_{k=1}^{2m}(q+(e^{ik\theta_{m}}+e^{-ik\theta_{m}}+2)^{n})^{2}.
\end{eqnarray*}
From the definition of $K_{n}$, we can list all roots of this polynomial
in pairs $(x_{l},x_{l}^{-1})$ so that
\[
e^{il\phi_{n}+i\pi/n}\sqrt[n]{|q|}=x_{l}+x_{l}^{-1}+2
\]
where $\phi_{n}=2\pi/n$ and $1\le l\le n$. Thus 
\begin{eqnarray*}
q+(e^{ik\theta_{m}}+e^{-ik\theta_{m}}+2)^{n} & = & (e^{ik\theta_{m}}+e^{-ik\theta_{m}}+2)^{n}-(e^{i\pi/n}\sqrt[n]{|q|})^{n}\\
 & = & \prod_{l=1}^{n}((e^{ik\theta_{m}}+e^{-ik\theta_{m}}+2)-e^{il\phi_{n}+i\pi/n}\sqrt[n]{|q|})\\
 & = & \prod_{l=1}^{n}((e^{ik\theta_{m}}+e^{-ik\theta_{m}}+2)-(x_{l}+x_{l}^{-1}+2))\\
 & = & \prod_{l=1}^{n}\frac{(x_{l}-e^{ik\theta_{m}})(x_{l}^{-1}-e^{ik\theta_{m}})}{e^{ik\theta_{m}}}.
\end{eqnarray*}
Since $\prod_{1\le k\le2m}e^{ik\theta_{m}}=1$, the proof follows.

\medskip{}
It follows that all roots of $\Delta_{x}(K_{n}f_{m})$ stay in the
range $[-2^{n},0]$. This is a special case of Stolarsky's proposition
stated in the introduction section.

\paragraph*{Corollary.}

For any $m$, the $H_{3m+1}^{(n)}(q)$ polynomials are divisible by
$(q+1)$. \\
Proof. Consider 
\[
H_{3m+1}^{(n)}(q)=\prod_{k=1}^{3m+1}(q+(2\cos k\theta_{3m+1}+2)^{n})
\]
 where 
\[
\theta_{3m+1}=\frac{2\pi}{6m+3}\,\,.
\]
When $k=2m+1$ we have $(2\cos k\theta_{3m+1}+2)^{n}=1$. Thus $(q+1)\,|\,H_{3m+1}^{(n)}(q)$.

\paragraph{Corollary.}

For any $n$, 
\[
H_{m}^{(n)}(q)\,|\,H_{3m+1}^{(n)}(q).
\]
\\
Proof. This is clear since the terms $k=3$, $6$, $9$, $\ldots$,
$3m$ in the product for $H_{3m+1}^{(n)}(q)$ give $H_{m}^{(n)}(q)$.

\section{Generating function for $H_{m}^{(1)}(q)$}

It is not hard to show that $H_{m}^{(1)}(q)$ has a generating function
similar to that of the closely related Chebyshev polynomials. We give
the details both for the sake of completeness, and because they are
useful in the rather harder analysis required for $H_{m}^{(2)}(q)$.

\paragraph{Proposition 3.1. }

The polynomials $H_{m}^{(1)}(x)$ satisfy
\[
\frac{1-t}{(1-t)^{2}-xt}=1+\sum_{m=1}^{\infty}H_{m}^{(1)}(x)t^{m}\,.
\]
Proof. Recall the generating function definition of Chebyshev T-polynomial
$T_{n}(x)$:
\[
\frac{1-xt}{1+t^{2}-2tx}=1+\sum_{m=1}^{\infty}T_{m}(x)t^{m}.
\]
By replacing $x$ by $(x+2)/2$, multiplying both sides by 2 and subtracting
1 from each side we obtain
\[
\frac{1-t}{(1-t)^{2}-xt}=\frac{1}{1+t}\left(1+\sum_{m=1}^{\infty}2T_{m}\left(\frac{x+2}{2}\right)t^{m}\right)\,.
\]
So it remains to prove that $H_{m}^{(1)}(x)$ is twice an alternating
sum of Chebyshev polynomials $T_{m}((x+2)/2)$ plus or minus 1 depending
on the parity of $m$, i.e.
\begin{equation}
H_{m}^{(1)}(x)=2T_{m}\left(\frac{x+2}{2}\right)-2T_{m-1}\left(\frac{x+2}{2}\right)+2T_{m-2}\left(\frac{x+2}{2}\right)-\cdots\mp2T_{1}\left(\frac{x+2}{2}\right)\pm1.
\end{equation}
By writing $2\cos m\theta=z^{m}+z^{-m}$ where $z=e^{i\theta}$, we
obtain
\[
S_{m}(z):=2\cos m\theta-2\cos(m-1)\theta+2\cos(m-2)\theta-\cdots\mp2\cos\theta\pm1=\frac{z^{2m+1}+1}{z^{m}(z+1)}\,.
\]
$S_{m}(z)$ has roots $e^{-ik\theta_{m}}$, so the twice alternating
sum of Chebyshev polynomial plus or minus 1 have roots $x$ such that
\[
\frac{x+2}{2}=-\cos k\theta_{m}
\]
or $x=-2\cos k\theta_{m}-2.$ These roots correspond to roots of $H_{m}^{(1)}(x)$.
This completes the proof by noting that both sides of (2) are monic
polynomials in $x$. 

\medskip{}

Proposition 3.1 leads to an explicit formula for the coefficients
of $H_{m}^{(1)}(x)$. This formula is not new, but will be useful
in finding the generating function for $H_{m}^{(2)}(x)$ which will
be discussed in section 3.

\paragraph*{Proposition 3.2.}

The polynomial $H_{m}^{(1)}(x)$ is given by the formula
\[
H_{m}^{(1)}(x)=\sum_{k=0}^{m}\left(\begin{array}{c}
m+k\\
m-k
\end{array}\right)x^{k}.
\]
Proof. Let $a_{k,m}$ be the coefficients of the polynomial $H_{m}^{(1)}(x)$.
Then by interchanging summation, we have the following identity:
\begin{eqnarray*}
\sum_{m=0}^{\infty}H_{m}^{(1)}(x)t^{m} & = & \sum_{m=0}^{\infty}\sum_{k=0}^{m}a_{k,m}x^{k}t^{m}\\
 & = & \sum_{k=0}^{\infty}x^{k}t^{k}\sum_{m=0}^{\infty}a_{k,m+k}t^{m}.
\end{eqnarray*}
 On the other hand, by expanding the generating function for $H_{m}^{(1)}(x)$
in terms of $x$, we obtain

\begin{eqnarray*}
\frac{1-t}{(1-t)^{2}-xt} & = & \frac{1}{1-t}\frac{1}{1-xt/(1-t)^{2}}\\
 & = & \sum_{k=0}^{\infty}x^{k}t^{k}/(1-t)^{2k+1}\\
 & = & \sum_{k=0}^{\infty}x^{k}t^{k}\sum_{m=0}^{\infty}\left(\begin{array}{c}
m+2k\\
m
\end{array}\right)t^{m}.
\end{eqnarray*}
 By equating the two double summations, we find that
\[
a_{k,m+k}=\left(\begin{array}{c}
m+2k\\
m
\end{array}\right),
\]
and the proposition follows.

\section{Generating function for $H_{m}^{(2)}(x)$}

Stolarsky \cite{stolarsky} conjectures that the generating function
of $H_{m}^{(2)}(x)$ will have the form 
\[
\frac{(1-t)^{3}}{(1-t)^{4}-xt(t+1)^{2}}.
\]
Our main approach to prove this formula is to express $H_{m}^{(2)}(x)$
in terms of $H_{m}^{(1)}(x)$ and apply proposition 3.2 in the previous
section. In particular, the connection between $H_{m}^{(2)}(x)$ and
$H_{m}^{(1)}(x)$ is given by 
\[
H_{m}^{(2)}(-x^{2})=H_{m}^{(1)}(x)H_{m}^{(1)}(-x).
\]
This equation easily follows from the definition of $H_{m}^{(n)}(x)$.
For the rest of this section, we use the following identity:

\[
\sum_{i=0}^{2k}(-1)^{i}\left(\begin{array}{c}
m+k+i\\
m+k-i
\end{array}\right)\left(\begin{array}{c}
m+3k-i\\
m-k+i
\end{array}\right)=(-1)^{k}\sum_{i=0}^{m}\left(\begin{array}{c}
2k\\
m-i
\end{array}\right)\left(\begin{array}{c}
4k+i\\
i
\end{array}\right).\tag{\ensuremath{\ast}}
\]
We will provide a proof for this identity in section 5.

\paragraph*{Proposition 4.1.}

The polynomials $H_{m}^{(2)}(x)$ satisfy
\[
\frac{(1-t)^{3}}{(1-t)^{4}-xt(t+1)^{2}}=1+\sum_{m=1}^{\infty}H_{m}^{(2)}(x)t^{m}.
\]
Proof. Recall that $H_{m}^{(2)}(-x^{2})=H_{m}^{(1)}(x)H_{m}^{(1)}(-x)$.
So $H_{m}^{(2)}(-x^{2})$ is an even polynomial of degree $2m$, whose
$k^{th}$-coefficient is given by
\[
\sum_{i=0}^{k}(-1)^{i}\left(\begin{array}{c}
m+i\\
m-i
\end{array}\right)\left(\begin{array}{c}
m+k-i\\
m-k+i
\end{array}\right).
\]
By interchanging summations as in the proof of Proposition 3.2, we
obtain
\[
\sum_{m=0}^{\infty}H_{m}^{(2)}(-x^{2})t^{m}=\sum_{k=0}^{\infty}x^{2k}t^{k}\sum_{m=0}^{\infty}\sum_{i=0}^{2k}(-1)^{i}\left(\begin{array}{c}
m+k+i\\
m+k-i
\end{array}\right)\left(\begin{array}{c}
m+3k-i\\
m-k+i
\end{array}\right)t^{m}.
\]
Upon expanding the function 
\[
\frac{(1-t)^{3}}{(1-t)^{4}+x^{2}t(t+1)^{2}}
\]
in terms of $x$ first and then in terms of $t$ one obtains 
\begin{eqnarray*}
\frac{(1-t)^{3}}{(1-t)^{4}+x^{2}t(t+1)^{2}} & = & \sum_{k=0}^{\infty}(-1)^{k}x^{2k}t^{k}(t+1)^{2k}/(1-t)^{4k+1}\\
 & = & \sum_{k=0}^{\infty}(-1)^{k}x^{2k}t^{k}\sum_{i=0}^{2k}\left(\begin{array}{c}
2k\\
i
\end{array}\right)t^{i}\sum_{j=0}^{\infty}\left(\begin{array}{c}
4k+j\\
j
\end{array}\right)t^{j}\\
 & = & \sum_{k=0}^{\infty}(-1)^{k}x^{2k}t^{k}\sum_{m=0}^{\infty}\sum_{i=0}^{m}\left(\begin{array}{c}
2k\\
m-i
\end{array}\right)\left(\begin{array}{c}
4k+i\\
i
\end{array}\right)t^{m}.
\end{eqnarray*}
 The proposition follows from the identity ($\ast$) above.

\medskip{}

\textit{Remark. }It is interesting to note that that since $H_{m}^{(2)}(-x^{2})$
is an even function, the following identity holds for any odd integer
$k$:
\[
\sum_{i=0}^{2m}(-1)^{i}\left(\begin{array}{c}
m+i\\
m-i
\end{array}\right)\left(\begin{array}{c}
m+k-i\\
m-k+i
\end{array}\right)=0.
\]

\section{A hypergeometric identity via Euler's contiguous relation and the
Wilf-Zeilberger algorithm }

In this section we will derive a proof for the identity
\[
\sum_{i=0}^{2k}(-1)^{i}\left(\begin{array}{c}
m+k+i\\
m+k-i
\end{array}\right)\left(\begin{array}{c}
m+3k-i\\
m-k+i
\end{array}\right)=(-1)^{k}\sum_{i=0}^{m}\left(\begin{array}{c}
2k\\
m-i
\end{array}\right)\left(\begin{array}{c}
4k+i\\
i
\end{array}\right).\tag{\ensuremath{\ast}}
\]
The method of proving this identity is similar to that of Vid$\bar{\mbox{u}}$nas
\cite{vidunas}. We first express the right hand side in terms of
hypergeometric $_{2}F_{1}$ function.

\paragraph*{Proposition 5.1.}

The following equations hold:
\[
\sum_{i=0}^{m}\left(\begin{array}{c}
2k\\
m-i
\end{array}\right)\left(\begin{array}{c}
4k+i\\
i
\end{array}\right)=\left(\begin{array}{c}
2k\\
m
\end{array}\right){}_{2}F_{1}\left(\begin{array}{c}
4k+1,-m\\
2k-m+1
\end{array};-1\right)
\]
if $2k\ge m$ and
\[
\sum_{i=0}^{m}\left(\begin{array}{c}
2k\\
m-i
\end{array}\right)\left(\begin{array}{c}
4k+i\\
i
\end{array}\right)=\left(\begin{array}{c}
m+2k\\
m-2k
\end{array}\right){}_{2}F_{1}\left(\begin{array}{c}
-2k,2k+m+1\\
m-2k+1
\end{array};-1\right)
\]
if $2k\le m$.

Proof. The first equation follows directly from the fact that
\[
\left(\begin{array}{c}
4k+i\\
i
\end{array}\right)=\frac{(4k+1)_{i}}{i!}
\]
and
\[
\left(\begin{array}{c}
2k\\
m-i
\end{array}\right)=\frac{(-1)^{i}(2k)!(-m)_{i}}{m!(2k-m)!(2k-m+1)!}.
\]
For the second equation, we note that the summand on the left side
equals 0 when $i<m-2k$. Thus the left side equals
\[
\sum_{i=m-2k}^{m}\left(\begin{array}{c}
2k\\
m-i
\end{array}\right)\left(\begin{array}{c}
4k+i\\
4k
\end{array}\right)=\sum_{i=0}^{2k}\left(\begin{array}{c}
2k\\
i
\end{array}\right)\left(\begin{array}{c}
4k+m-i\\
4k
\end{array}\right)=\sum_{i=0}^{2k}\left(\begin{array}{c}
2k\\
i
\end{array}\right)\left(\begin{array}{c}
m+2k+i\\
4k
\end{array}\right).
\]
To complete the proof note that
\[
\left(\begin{array}{c}
2k\\
i
\end{array}\right)=\frac{(-2k)_{i}(-1)^{i}}{i!}
\]
and
\[
\left(\begin{array}{c}
2k+m+i\\
4k
\end{array}\right)=\frac{(2k+m)!(2k+m+1)_{i}}{(4k)!(m-2k)!(m-2k+1)_{i}}.
\]
\medskip{}

To continue the proof of ($\ast$), we recall Euler's contiguous relation:
\begin{eqnarray*}
_{2}F_{1}\left(\begin{array}{c}
a,b\\
c
\end{array};u\right) & = & \frac{c+(a-b+1)u}{c}{}_{2}F_{1}\left(\begin{array}{c}
a+1,b\\
c+1
\end{array};u\right)\\
 &  & -\frac{(a+1)(c-b+1)u}{c(c+1)}{}_{2}F_{1}\left(\begin{array}{c}
a+2,b\\
c+2
\end{array};u\right).
\end{eqnarray*}
Applying this relation to 
\[
_{2}F_{1}\left(\begin{array}{c}
4k+1,-m\\
2k-m+1
\end{array};-1\right)
\]
with $a=-m$, $b=4k+1$ and $c=2k-m+1$, we obtain the identity
\begin{eqnarray*}
_{2}F_{1}\left(\begin{array}{c}
4k+1,-m\\
2k-m+1
\end{array};-1\right) & = & \frac{6k+1}{2k-m+1}{}_{2}F_{1}\left(\begin{array}{c}
4k+1,-m+1\\
2k-m+2
\end{array};-1\right)\\
 &  & +\frac{(m-1)(m-1+2k)}{(2k-m+1)(2k-m+2)}{}_{2}F_{1}\left(\begin{array}{c}
4k+1,-m+2\\
2k-m+3
\end{array};-1\right).
\end{eqnarray*}
Also applying the same contiguous relation to 
\[
_{2}F_{1}\left(\begin{array}{c}
-2k,2k+m+1\\
m-2k+1
\end{array}\right)
\]
with $a=2k+m+1,$ $b=-2k$ and $c=m-2k+1$, we obtain
\begin{eqnarray*}
_{2}F_{1}\left(\begin{array}{c}
-2k,2k+m+1\\
m-2k+1
\end{array};-1\right) & = & -\frac{6k+1}{m-2k+1}{}_{2}F_{1}\left(\begin{array}{c}
-2k,2k+m+2\\
m-2k+2
\end{array};-1\right)\\
 &  & +\frac{(m+2)(2k+m+2)}{(m-2k+1)(m-2k+2)}{}_{2}F_{1}\left(\begin{array}{c}
-2k,2k+m+3\\
m-2k+3
\end{array};-1\right).
\end{eqnarray*}
Fix $k$ and denote by $S_{m}$ the right side of $(\ast)$. From
the two identities above and Proposition 5.1, we have the recursive
relation
\[
(m+2)S_{m+2}=(6k+1)S_{m+1}+(m+1+2k)S_{m},
\]
where $S_{0}=(-1)^{k}$ and $S_{1}=(-1)^{k}(6k+1)$.

Let $T_{m}$ be the left side of $(\ast)$. It is easy to check that
$T_{0}=(-1)^{k}$ and
\begin{eqnarray*}
T_{1} & = & (-1)^{k-1}\left(\begin{array}{c}
2k\\
2
\end{array}\right)+(-1)^{k}(2k+1)^{2}+(-1)^{k+1}\left(\begin{array}{c}
2k\\
2
\end{array}\right)\\
 & = & (-1)^{k}(6k+1).
\end{eqnarray*}
 So it suffices to show that the sequence $T_{m}$ also satisfies
the same recursive relation, i.e.
\[
(m+2)T_{m+2}-(6k+1)T_{m+1}-(m+1+2k)T_{m}=0.
\]
By definition of $T_{m}$, the left hand side of the relation above
is a finite summation of several terms. We denote by $f(m,i)$ its
summand. Also let $F(m,i)=f(m,i)/(2k+m+1)!$. Following the WZ-algorithm
we can define the certificate function
\begin{eqnarray*}
R(m,i) & := & \frac{i(2i-1)(i-3k-m-1)}{2(-3+i-k-m)(2+2k+m)}\frac{R_{1}(m,i)}{R_{2}(m,i)}
\end{eqnarray*}
where
\begin{eqnarray*}
R_{1}(m,i) & = & -10-i+i^{2}-5k-4ik+2i^{2}k+7k^{2}-4ik^{2}\\
 &  & +2k^{3}-13m-im+i^{2}m-5km\\
 &  & -2ikm+k^{2}m-6m^{2}-2km^{2}-m^{3}
\end{eqnarray*}
and
\begin{eqnarray*}
R_{2}(m,i) & = & -5i^{2}+2i^{4}+2k+10ik-3i^{2}k-8i^{3}k\\
 &  & -5k^{2}+6ik^{2}+6i^{2}k^{2}-11k^{3}+4ik^{3}-4k^{4}\\
 &  & -6i^{2}m+3km+12ikm-4i^{2}km-10k^{2}m+8ik^{2}m\\
 &  & -8k^{3}m-2i^{2}m^{2}+km^{2}+4ikm^{2}-4k^{2}m^{2}.
\end{eqnarray*}
 With computer algebra one can check that 
\[
F(m+1,i)-F(m,i)=G(m,i+1)-G(m,i)
\]
where $G(m,i)=R(m,i)F(m,i)$. The equation above can be checked easily
 by Mathematica using the FactorialSimplify function in the aisb.m
package (see \cite{petkovsek}). We note that $G(m,0)=0$ since $R(m,0)=0$
by definition of $R(m,i)$. Also $G(m,2k+1)=0$ since $F(m,2k+1)=0$.
Therefore 
\[
\sum_{i=0}^{2k}F(m,i)
\]
 is a constant and it is easy to check that this constant is 0 by
the initial condition. Thus 
\[
\sum_{i=0}^{2k}f(m,i)=0,
\]
and $T_{m}$ satisfies the recursive relation. \\
\\


\begin{thebibliography}{10}
\bibitem[1]{aar} G. E. Andrews, R. Askey, and R. Roy, $\textit{Special Functions}$,
Cambridge University Press, Cambridge, 1999.

\bibitem[2]{apostol} T. M. Apostol, The resultants of the cyclotomic
polynomials $F_{m}(ax)$ and $F_{n}(bx)$, Math. Comp. 29 (1975),
1-6.

\bibitem[3]{dilcherstolarsky} K. Dilcher and K. B. Stolarsky, Resultants
and Discriminants of Chebyshev and related polynomials, Transactions
of the Amer. Math. Soc. 357 (2004), 965-981.

\bibitem[4]{gkz} I. M. Gelfand, M. M. Kapranov, and A. V. Zelevinsky,
$\mathit{\textit{Discriminants, Resultants, and}}$ $\textit{Multidimensional Determinants}$,
Birkhuser Boston, Boston, 1994.

\bibitem[5]{gishe-ismail} J. Gishe and M. E. H. Ismail, Resultants
of Chebyshev Polynomials, Z. Anal. Anwend. 27 (2008), no. 4, 499-508.

\bibitem[6]{petkovsek} M. Petkovsek, Computer algebra package aisb.m
for Mathematica. http://www.fmf.uni-lj.si/\textasciitilde{}petkovsek/software.html.

\bibitem[7]{pwz} M. Petkovsek, H. S. Wilf, and D. Zeilberger, $\textit{A=B}$.
A K Peters, Ltd, 1996.

\bibitem[8]{roberts} D. P. Roberts, Discriminants of some Painlev$\acute{\mbox{e}}$
polynomials, Number theory for the millennium, III, A K Peters, Natick,
2002, pp. 205-221.

\bibitem[9]{stolarsky} K. Stolarsky, Discriminants and divisibility
for Chebyshev-like polynomials, Number theory for the millennium,
III (Urbana, IL, 2000), 243--252, A K Peters, Natick, MA, 2002.

\bibitem[10]{vidunas} R. Vid$\bar{\mbox{u}}$nas, A generalization
of Kummer's identity, Conference on Special Functions (Tempe, AZ,
2000). Rocky Mountain J. Math. 32 (2002), no. 2, 919--936. \end{thebibliography}
\end{document}